\documentclass[12pt]{amsart}
\pdfoutput=1
\usepackage{amsmath}
\usepackage{amsfonts}
\usepackage{amsthm}
\usepackage{graphicx}
\usepackage{color}
\usepackage[margin=1.01in]{geometry}
\usepackage{subfig}
\usepackage{float}
\usepackage{tikz}
\usepackage{enumerate}
\usepackage[titletoc]{appendix}
\usepackage{microtype}
\theoremstyle{definition}

\newtheorem{defn}{Definition}[section]

\theoremstyle{Theorem}
\newtheorem{thm}{Theorem}[section]
\newtheorem{lemma}[thm]{Lemma}
\theoremstyle{remark}

\newcommand{\e}{\varepsilon}

\newcommand{\ie}{\emph{i.e.}}

\newcommand{\pder}[3]{\frac{\partial^{#3} #1}{\partial #2^{#3}}}
\newcommand{\mL}{\mathcal{L}}

\usepackage{pgf}

\begin{document}
\title{Spectral Stability of Travelling Wave Solutions in a Keller-Segel Model}


\author{P.N. Davis$^\dagger$}
  \thanks{$^\dagger$School of Mathematical Sciences, Queensland University of Technology, Brisbane, Australia}
\author{P. van Heijster$^\dagger$}
 \author{R. Marangell$^\ddagger \ ^*$}
 \thanks{$^\ddagger$School of Mathematics and Statistics, University of Sydney, Sydney, Australia}
 \thanks{ $^*$ Corresponding author: robert.marangell@sydney.edu.au}
\begin{abstract}
We investigate the point spectrum associated with travelling wave solutions in a Keller-Segel model for bacterial chemotaxis with small diffusivity of the chemoattractant, a logarithmic chemosensitivity function and a constant, sublinear or linear consumption rate. We show that, for constant or sublinear consumption, there is an eigenvalue at the origin of order two. This is associated with the translation invariance of the model and the existence of a continuous family of solutions with varying wave speed. These point spectrum results, in conjunction with previous results in the literature, imply that in these cases the travelling wave solutions are absolute unstable if the chemotactic coefficient is above a certain critical value, while they are transiently unstable otherwise.
 
\end{abstract}
\maketitle
\section{Introduction}
The aim of this manuscript is to complete the spectral results for travelling wave solutions to the nondimensionalised Keller-Segel model for bacterial chemotaxis with zero-growth (or decay) of the bacterial population $w$ and a logarithmic chemosensitivity function
\begin{equation}
\label{EQ:KSmain}
\begin{aligned}
u_{{t}}&={\e} u_{{x}{x}}- w u^m,\\
w_{{t}}&= w_{{x}{x}}-{\beta}\left( \Phi_x(u) w\right)_{{x}}{\,,} 
\quad {\textnormal{with}} \quad \Phi(u) = \ln{u}\,.
\end{aligned}
\end{equation}
Here, $0 \leq \e \ll 1$ represents the small diffusivity of the chemoattractant $u$, 
$(x,t)\in\mathbb{R}\times\mathbb{R}^+$, $0\leq m\leq 1$, and $\beta>1-m$. This nondimensionalised
form \eqref{EQ:KSmain} was introduced in \cite{nagai1991traveling} for $m=0$ and used in \cite{davis2017absolute} for $0\leq m \leq 1$.

The original Keller-Segel model was proposed by Keller and Segel in the 1970's \cite{keller1971model, keller1971traveling}. Much of the existing literature has focused on the minimal Keller-Segel model with linear chemosensitivity functions such as $\Phi(u)=u$ 
\cite{horstmann20031970,horstmann2005boundedness,kang2007stability} and on Keller-Segel models with singular chemosensitivity functions such as $\Phi(u)=\ln u$ \cite{ebihara1992singular,feltham2000travelling,keller1975necessary,keller1971traveling, schwetlick2003traveling}. For instance, it was shown in \cite{keller1971traveling, schwetlick2003traveling} that, in the absence of growth term for the bacterial population $w$, the chemosensitivity function must be singular for travelling wave solutions to exist. Hence our choice of a logarithmic chemosensitivity function in \eqref{EQ:KSmain}.

 The existence of travelling wave solutions that travel with constant speed $c$ has been shown under the above parameter conditions for all $\e \geq 0$ \cite{feltham2000travelling, nagai1991traveling, wang2013mathematics}, \ie\ travelling wave solutions exist for $\e$ not necessarily small. However, the stability results of \cite{davis2017absolute} related to these travelling wave solutions only hold for small $\e$ and therefore we consider $0\leq \e\ll1$.
Linearising about a travelling wave solution in its associated moving frame results in a linearised operator $\mL_\e^m$ whose spectrum determines the spectral stability of the travelling wave solution.
The spectrum of the linearised operator falls naturally into two parts, the essential spectrum $\sigma_{\rm ess}$ and the point spectrum $\sigma_{\rm pt}$, see \cite{kapitula2013spectral,sandstede2002stability} and references therein.
 In \cite{davis2017absolute}, it was shown that the essential spectrum of $\mL_\e^m$ is contained in the open left half plane in a two-sided exponentially weighted function space when $\beta<\beta_{\rm crit}^m(\e)$, where $\beta_{\rm crit}^m(\e)$ was found to be, to leading order in $\e$, the unique root $\beta$ greater than $1-m$ of a $10th$-order polynomial, see upcoming Theorem \ref{TH:MAIN1}. In \cite{harley2015numerical}, \eqref{EQ:KSmain} was studied with $m=0=\e$ and it was shown, via a numerical Evans function computation, that the point spectrum of $\mL_0^0$ contains no eigenvalues with positive real part for complex values with norm up to $\mathcal{O}(10^9)$. In addition, it was shown that the origin is a second order root of the Evans function and hence is an eigenvalue with algebraic multiplicity two.

In this manuscript we prove that the origin persists as an element of the point spectrum with algebraic multiplicity two for travelling wave solutions to \eqref{EQ:KSmain} with $0\leq m<1$ and $0\leq \e\ll 1$. This is done by explicitly solving the associated (generalised) eigenvalue problem and analysing the asymptotic behaviour of the (generalised) eigenfunctions. In \S\ref{SS:setup}, we formulate the eigenvalue problem and outline the relevant spectral theory and existing stability results for travelling wave solutions to \eqref{EQ:KSmain}. In \S \ref{SS:point_spectrum}, we study the eigenvalue problem by first computing the eigenfunction of the origin associated with the translation invariance of \eqref{EQ:KSmain}. Next, we compute a generalised eigenfunction, which is due to the existence of a family of travelling wave solutions under varying the wave speed $c$. By analysing the asymptotic behaviour of the (generalised) eigenfunctions we prove that they are contained in a range of exponentially weighted function spaces for which the essential spectrum is stable. We conclude the spectral stability of travelling wave solutions to \eqref{EQ:KSmain} in these function spaces. In \S \ref{SS:linear_m} we extend the analysis to the $m=1$ case. In this case we can not conclude spectral stability as there is no spectral gap in any exponentially weighted function space. 
However, nonlinear (in)stability results have been obtained in certain cases for Keller-Segel models with linear consumption rate \cite{li2014stability,li2012steadily,meyries2011local}. We include the $m = 1$ case for completeness and to highlight how and why the spectral gap vanishes in the $m \to 1$ limit.
We conclude the manuscript with a summary and discussion of future work.

\section{Set-up of the Eigenvalue Problem}\label{SS:setup}

After the change of variables $z=x-ct$, where $c>0$ is the constant, finite wave speed, \eqref{EQ:KSmain} becomes
\begin{equation}
\label{EQ:KStw}
\begin{split}
u_t&=\e u_{zz}+cu_z- w u^m,\\
w_t&=w_{zz}+cw_z-\beta\left(\frac{wu_z}{u}\right)_z \,.
\end{split}
\end{equation}
For $0\leq m<1$, travelling wave solutions are stationary solutions to \eqref{EQ:KStw}, \ie\ $(u(z,t),w(z,t)) = (u^m_\e(z),w^m_\e(z))$, that asymptote to
\begin{equation}\label{EQ:limits1}
\lim_{z\to-\infty} (u^m_\e(z),w^m_\e(z))=(0,0),\quad\lim_{z\to\infty}(u^m_\e(z),w^m_\e(z))=(u_r,0)\,.
\end{equation}
Without loss of generality, $u_r$ can be scaled to one \cite{nagai1991traveling, wang2013mathematics}. 
In other words, $u^m_\e(z)$ is a wavefront and $w^m_\e(z)$ is a pulse and they satisfy
\begin{equation}
\label{EQ:KStw2}
\begin{split}
0 &= \e u_{zz}+cu_z- w u^m,\\
0 &= w_{zz}+cw_z-\beta\left(\frac{wu_z}{u}\right)_z .
\end{split}
\end{equation}
The existence of these travelling wave solutions $(u^m_\e,w^m_\e)$ was proven for $\e\geq0$ in \cite{feltham2000travelling,nagai1991traveling,wang2013mathematics}. For  
$\e$ small, the profiles are to leading order given by
\begin{equation}
\label{tw_profiles_0}
\begin{split}
u_0^m(z)=\left( 1+\frac{(\beta+m-1)}{c^2} e^{-cz} \right)^{-1/(\beta+m-1)}\,,\qquad
w_0^m(z)=e^{-cz}\left(u^m_0(z)\right)^{\beta}\,,
\end{split}
\end{equation}
where we, without loss of generality, centred the travelling wave solution around $z=0$. Observe that the relationship between $w_0^m$ and $u_0^m$ also hold for $\e \neq 0$. That is, $w_\e^m(z)=e^{-cz}\left(u^m_\e(z)\right)^{\beta}$ \cite{feltham2000travelling}.

To determine the spectral stability of these travelling wave solutions $(u^m_\e,w^m_\e)$, we substitute $u(z,t)=u^m_\e(z)+p(z)e^{\lambda t}$ and $w(z,t)=w^m_\e(z)+q(z)e^{\lambda t}$ into \eqref{EQ:KStw}. Here, $p$ and $q$ are perturbations in $\mathbb{H}^1(\mathbb{R})$, the Sobolev space of once (weakly) differentiable functions such that both the function and its first (weak) derivative (in $z$) are in $\mathbb{L}^2(\mathbb{R})$, \ie\ both are square integrable. By considering only leading order terms for $p$ and $q$, we obtain the linear operator $\mL_\e^m:\mathbb{H}^1(\mathbb{R})\times \mathbb{H}^1(\mathbb{R})\to \mathbb{H}^1(\mathbb{R})\times \mathbb{H}^1(\mathbb{R})$ and the eigenvalue problem given by 
 \begin{align}&\mL_\e^m \begin{pmatrix}
 p\\q
 \end{pmatrix}=\lambda\begin{pmatrix}
 p\\q
 \end{pmatrix},\quad\quad
\mathcal{L}_\e^m:=\begin{pmatrix}\mathcal{L}_{11} &\mathcal{L}_{12}\\ \mathcal{L}_{21} &\mathcal{L}_{22}.
\end{pmatrix}.\label{EQ:Lop_sublinear_m}
\end{align}
The entries of $\mathcal{L}_\e^m$ are
\begin{equation}
\label{EQ:Lop_sublinear_m_entries}
\begin{aligned}
\mathcal{L}_{11}&:=\e\pder{}{z}{2}+c\pder{}{z}{}- m w u^{m-1},\\
\mathcal{L}_{12}&:= -u^m,\\
\mathcal{L}_{21}&:=\beta\left(\frac{w_z u_z }{u^2}+\frac{wu_{zz} }{u^2}-\frac{2wu_{z}^2}{u^3}\right)+\beta\left(\frac{2wu_{z}}{u^2}-\frac{w_z}{u}\right)\pder{}{z}{}-\frac{\beta w}{u}\pder{}{z}{2},\\
\mathcal{L}_{22}&:=\beta\left(\frac{u_{z}^2}{u^2}-\frac{u_{zz}}{u}\right)+\left(c-\frac{\beta u_{z}}{u}\right)\pder{}{z}{}+\pder{}{z}{2},
\end{aligned}
\end{equation}
where we have dropped the $\e$ subscripts and $m$ superscripts from $(u_\e^m,w_\e^m)$ for convenience.

The spectrum of an operator consists of values of $\lambda\in\mathbb{C}$ such that the inverse of the eigenvalue operator $\mL -\lambda I$ does not exist or is unbounded. More precisely, 
\begin{defn}\label{defn:point_spectrum}(\cite{kapitula2013spectral} Definition 2.2.3)
For a closed linear operator $\mL:\mathcal{D}(\mL)\subset X\to X$, where $X$ is a Banach space and $\mathcal{D}(\mathcal{L})$ is dense in $X$, the spectrum is decomposed into two sets:
\begin{enumerate}[(a)]
\item The essential spectrum $\sigma_{\rm ess}$ of the operator $\mL$ is the set of all $\lambda\in\mathbb{C}$ such that 
\begin{itemize}
\item $\mL-\lambda I$ is not Fredholm or,
\item $\mL-\lambda I$ is Fredholm with a non zero Fredholm index,
\end{itemize} 
where the Fredholm index of $\mL$ is $$\text{ind}(\mL)=\text{dim(ker(}\mL))-\text{codim(range(}\mL)).$$
\item The point spectrum of the operator $\mL$ is the set of all $\lambda\in\mathbb{C}$ such that the operator $\mL-\lambda I$ is Fredholm with index zero, but the operator is not invertible. That is,
$$\sigma_{\rm pt}=\left\lbrace \lambda\in\mathbb{C}:\text{ind}(\mL-\lambda I)=0, \text{ but }\mL-\lambda I\text{ is not invertible}\right\rbrace.$$
\end{enumerate}
\end{defn}
As in \cite{kapitula2013spectral}, we use the term \textit{eigenvalue} to refer to all values $\lambda$ in the spectrum of the operator whereas \textit{point spectrum} refers to isolated eigenvalues of finite multiplicity. Eigenvalues which have eigenfunctions that decay to zero as $z\to\pm\infty$ may be embedded in the essential spectrum and are thus not considered point spectrum.

Later on it will be more convenient to work with a fourth order ordinary differential equation (ODE) equivalent to the eigenvalue problem \eqref{EQ:Lop_sublinear_m} with terms consisting only of the travelling wave solution $u_\e^m$, the perturbation function $p$ and their derivatives. By rearranging \eqref{EQ:KStw2} and \eqref{EQ:Lop_sublinear_m}, we obtain, after some algebra, the following singular fourth order ODE

\begin{align}
\e p_{zzzz}- \mathcal{D}_{m,\e}p_{zzz}- \mathcal{C}_{m,\e}p_{zz}- \mathcal{B}_{m,\e} p_z - \mathcal{A}_{m,\e} p = 0, \label{EQ:Pop_m_nonzero_eps_nonzero}
\end{align}
with 
\begin{equation*}
\begin{split}
\mathcal{A}_{m,\e}:=&(\beta+m) \left(c^2+\lambda +\lambda  m\right)\frac{ u_z^2}{u^2}-2 c \lambda  m\frac{ u_z}{u}-c (\beta+m)\frac{ u_z
   u_{zz}}{u^2}-\lambda ^2\\
   &-\lambda  (\beta+m)\frac{ u_{zz}}{u}- c (\beta -2)(\beta+m)\frac{ u_z^3}{u^3} \\
   & + \e  \bigg(c (\beta+m)\frac{ u_z u_{zz}}{u^2}-(\beta -2) (\beta+m) \frac{u_z^2 u_{zz}}{u^3}
   -(\beta+m)\frac{
   u_{zz}^2}{u^2}-\lambda  m\frac{ u_{zz}}{u}\bigg), \\
\mathcal{B}_{m,\e}:=&2 c \lambda-\left(\beta  c^2+\lambda (\beta+2 m)\right)\frac{u_z }{u}+c (\beta -m-3) (\beta+m)\frac{ u_z^2}{u^2}+c (\beta+m)\frac{ u_{zz}}{u}\\
&+\e  \left((\beta -2) (\beta+m)\frac{ u_z u_{zz}}{u^2}-c   (\beta+m)\frac{ u_{zz}}{u}\right),\\
\mathcal{C}_{m,\e}:=&-c^2+c (2 (\beta+m) +m)\frac{
   u_z}{u}+\lambda+\e \bigg(\lambda-(m+1) (\beta+m)\frac{ u_z^2}{u^2}+c m\frac{ u_z}{u}\\
   & +2 (\beta+m) \frac{u_{zz}}{u}\bigg),\\
\mathcal{D}_{m,\e}:=&-c+\e \left((\beta+2 m) \frac{u_z}{u}-c\right).
\end{split}
\end{equation*}
Though the terms of \eqref{EQ:Pop_m_nonzero_eps_nonzero} appear singular as $u_\e^m\to0$ in the $z\to-\infty$ limit, they are in fact bounded \cite{wang2013mathematics}.

The essential spectrum of the operator $\mathcal{L}$ in \eqref{EQ:Lop_sublinear_m} depends on the asymptotic behaviour of the operator. In particular, it depends on the magnitude and signs of the spatial eigenvalues of the asymptotic states as $z\to\pm\infty$. These spatial eigenvalues are found as the roots of the following equations
\begin{equation}\label{EQ:characteristic_eqn}
\begin{split}
\e (\mu^+)^4- \mathcal{D}_{m,\e}^+(\mu^+)^3- \mathcal{C}_{m,\e}^+(\mu^+)^2- \mathcal{B}_{m,\e}^+(\mu^+)- \mathcal{A}_{m,\e}^+ &= 0, \\
\e (\mu^-)^4- \mathcal{D}_{m,\e}^-(\mu^-)^3- \mathcal{C}_{m,\e}^-(\mu^-)^2- \mathcal{B}_{m,\e}^-(\mu^-)- \mathcal{A}_{m,\e}^- &= 0 ,
\end{split}
\end{equation}
where $\mathcal{A}_{m,\e}^\pm$, $\mathcal{B}_{m,\e}^\pm$, $\mathcal{C}_{m,\e}^\pm$ and $\mathcal{D}_{m,\e}^\pm$ respectively denote the limits of $\mathcal{A}_{m,\e}$, $\mathcal{B}_{m,\e}$, $\mathcal{C}_{m,\e}$ and $\mathcal{D}_{m,\e}$ as $z\rightarrow\pm\infty$. Observe that the expressions \eqref{EQ:characteristic_eqn} are exactly the characteristic equations of \eqref{EQ:Pop_m_nonzero_eps_nonzero} in the limit $z\to\pm\infty$. The essential spectrum consists of values $\lambda\in\mathbb{C}$ such that any of the spatial eigenvalues $\mu^\pm$ are purely imaginary or the number of unstable spatial eigenvalues at $\pm\infty$ differ. This is equivalent to Definition 
\ref{defn:point_spectrum} \cite{kapitula2013spectral}. 

All travelling wave solutions to \eqref{EQ:KStw} have essential spectrum in the right half plane for perturbations in $\mathbb{H}^1(\mathbb{R})$ \cite{davis2017absolute,nagai1991traveling}. Thus, we follow the usual procedure outlined in \cite{kapitula2013spectral} and introduce 
the weighted function space $\mathbb{H}^1_\nu(\mathbb{R})$ defined by the norm
\begin{align}
\|p\|_{\mathbb{H}^1_\nu} =\|e^{\nu z}p\|_{\mathbb{H}^1}=\|\tilde{p}\|_{\mathbb{H}^1},\label{EQ:weighted_norm}
\end{align}
where $\tilde{p}:=e^{\nu z} p$. So, $p\in \mathbb{H}^1_\nu$ if and only if $\tilde{p}\in \mathbb{H}^1$. We define $\mathbb{L}^2_\nu$ similarly. 
It was shown in \cite{davis2017absolute} that a two-sided weight is required for the current problem. That is, 
\begin{equation}\label{EQ:two_sided_weight}\nu=\begin{cases} \nu_- \mbox{ if} & z \leq0, \\ \nu_+\mbox{ if} & z>0,\end{cases} \end{equation}
which forces the perturbation to decay exponentially in both directions. Using a weighted function space has the effect of shifting the essential spectrum. In particular, in the weighted function space we consider the spatial eigenvalues are $\mu^++\nu_+$ and $\mu^-+\nu_-$ as $z\to\pm\infty$ respectively. In other words, the weighted essential spectrum consists of values $\lambda\in\mathbb{C}$ such that any of the spatial eigenvalues $\mu^\pm+\nu_\pm$ are purely imaginary or the number of weighted unstable spatial eigenvalues at $\pm\infty$ differ.

  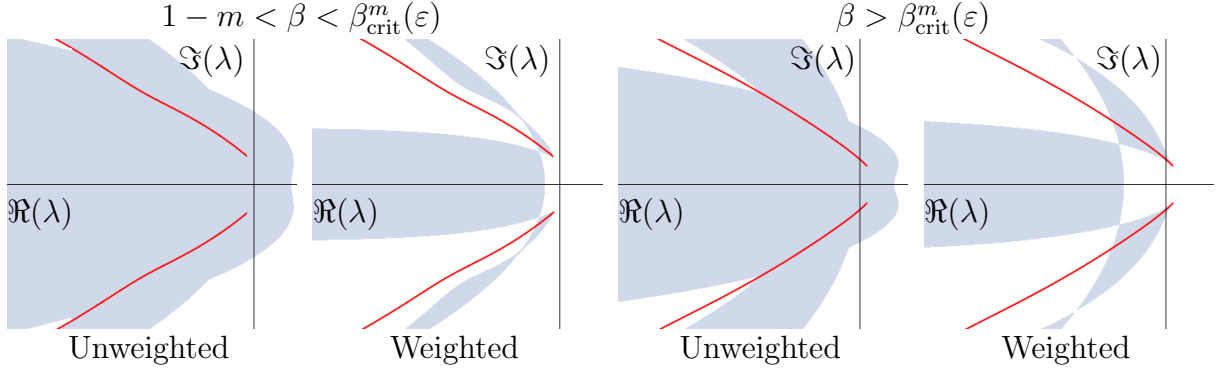
\begin{figure}[t]
\centering
\resizebox{\textwidth}{!}{
\begin{tikzpicture}
\node[anchor=south west] at (0,0){\includegraphics[width=.23\textwidth]{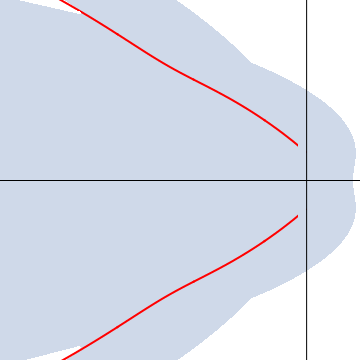}};
\node[anchor=north west] at (0,2){$\Re(\lambda)$};
\node[anchor=south west] at (4,0){\includegraphics[width=.23\textwidth]{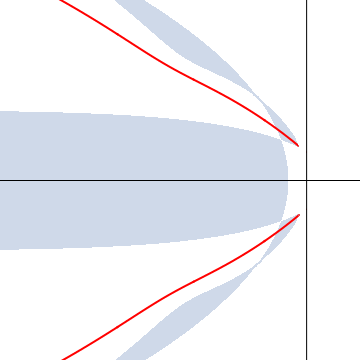}};
\node[anchor=north west] at (4,2){$\Re(\lambda)$};
\node[anchor=south west] at (8,0){\includegraphics[width=.23\textwidth]{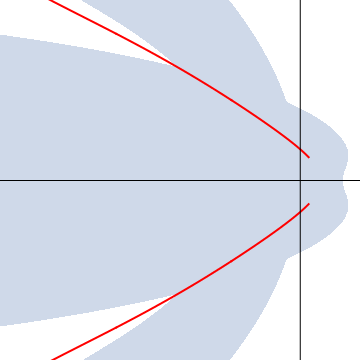}};
\node[anchor=north west] at (8,2){$\Re(\lambda)$};
\node[anchor=south west] at (12,0){\includegraphics[width=.23\textwidth]{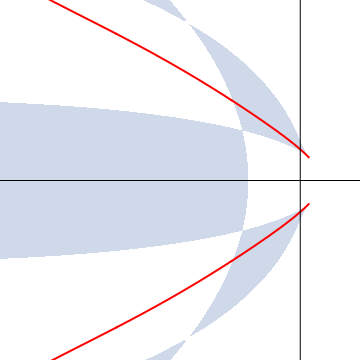}};
\node[anchor=north west] at (12,2){$\Re(\lambda)$};
\node[anchor=north east] at (3.4,4){$\Im(\lambda)$};
\node[anchor=north east] at (7.4,4){$\Im(\lambda)$};
\node[anchor=north east] at (11.4,4){$\Im(\lambda)$};
\node[anchor=north east] at (15.4,4){$\Im(\lambda)$};
\node[anchor=north] at (2,0.2){Unweighted};
\node[anchor=north] at (6,0.2){Weighted};
\node[anchor=north] at (10,0.2){Unweighted};
\node[anchor=north] at (14,0.2){Weighted};
\node[anchor=south] at (4,3.85){$1-m<\beta<\beta_{\rm crit}^m(\e)$};
\node[anchor=south] at (12,3.85){$\beta>\beta_{\rm crit}^m(\e)$};
\end{tikzpicture}}
\caption{A schematic of the typical weighted and unweighted essential spectra (blue regions) associated with travelling wave solutions $(u_\e^m,w_\e^m)$ of \eqref{EQ:KSmain}. The red curves indicate the subset of the absolute spectrum that determines how far the essential spectrum can be shifted by weighting the function space. For $1-m<\beta<\beta_{\rm crit}^m(\e)$ the absolute spectrum is contained in the left half plant and a two-sided weight exists such that the weighted essential spectrum is contained in the open left half plane.}
\label{FIG:schematic}
\end{figure} 
The following result is adapted from Theorem 2.6 of \cite{davis2017absolute}
\begin{thm}
\label{TH:MAIN1}
Assume that $c>0, 0 \leq m < 1$ and $\beta> 1-m$.
Let $\beta_{\rm crit}$ be the unique real root larger than one of 
\begin{equation}
\begin{aligned}
f(\beta) =\ &310 \beta^{10} - 3234 \beta^9 + 17112 \beta^8 - 49101 \beta^7 + 76180 \beta^6 - 58398 \beta^5 \\
& + 10056 \beta^4+ 15040 \beta^3 - 9680 \beta^2 + 1716 \beta -4. \label{EQ:10th_order_poly}
\end{aligned}
\end{equation}
Then, there exists an $\e_0>0$ such that for all $0\leq \e <\e_0$ there exists a range of two-sided weights \eqref{EQ:two_sided_weight} such that the weighted essential spectrum of  $\mathcal{L}_\e^m$ given in \eqref{EQ:Lop_sublinear_m} is fully contained in the left half plane for all $1-m<\beta<\beta_{\rm crit}^m(\e)$, with $\beta_{\rm crit}^m(\e)$ to leading order given by $\beta_{\rm crit}^m := \beta_{\rm crit}(1-m)$.
 For $\beta>\beta_{\rm crit}^m(\e)$ the essential spectrum of $\mathcal{L}_\e^m$ contains values in the right half plane for all possible weights and the travelling wave solutions of \eqref{EQ:KStw2} are thus absolutely unstable.
\end{thm}
This result was obtained by studying the so-called absolute spectrum of $\mL_\e^m$ which is not effected by weighting the function space (as it is related to the relative spacing between the real parts of the spatial eigenvalues). See Figure \ref{FIG:schematic} and \cite{davis2017absolute} for more detail.

All travelling wave solutions $ (u_\e^m(z),w_\e^m(z))$ are absolutely unstable for $\beta> \beta_{\rm crit}^m(\e)$, precluding the possibility of spectral stability. Thus we focus on the the parameter regime that is potentially transiently unstable, \ie\ $1-m<\beta<\beta_{\rm crit}^m(\e)$ in the remainder of this manuscript. To be able to conclude transient instability, \ie\ spectral stability in an exponentially weighted function space \cite{sandstede2002stability}, we must show that for the weights that shift the essential spectrum into the left half plane, there are no values $\lambda$ in the point spectrum with $\Re (\lambda)\geq0$ other than for $\lambda=0$. The location of the point spectrum does not change upon moving to a weighted space \cite{kapitula2013spectral}, however, it is necessary to show that the eigenfunctions associated with the point spectrum are contained in these weighted function spaces.

\section{Point Spectrum}\label{SS:point_spectrum}
Locating the point spectrum amounts to finding nontrivial solutions $(p,q)^T$ to \eqref{EQ:Lop_sublinear_m} that decay to zero as $z\to\pm\infty$ for some $\lambda\in\mathbb{C}\backslash \sigma_{\rm ess}$. While $\lambda=0$ is in the essential spectrum in the unweighted space for all $m$ it is not in the weighted function space for a range of weights when $0\leq m<1$ \cite{davis2017absolute,harley2015numerical}, see also Figure \ref{FIG:schematic}. It was shown in \cite{harley2015numerical} that $\lambda=0$ is a root with order two of the Evans function for travelling solutions to \eqref{EQ:KSmain} with $\e=m=0$. Thus, in these appropriately weighted function spaces $\lambda=0$ is an isolated eigenvalue associated to the invariances of the problem. Hence, $\lambda=0$ is part of the point spectrum. We show that the eigenvalue $\lambda=0$ persists with multiplicity two for $0\leq m<1$ by determining two linearly independent eigenfunctions that form the generalised eigenspace. In the singular limit $\e\to0$, these functions and their norms are explicitly computed. We also show these eigenfunctions persist for sufficiently small $\e>0$ and are contained in the weighted function spaces that have stable essential spectrum.\footnote{The formulation of the (generalised) eigenfunctions, eigenspace and analysis of the asymptotic behaviour of these functions remains valid when $\beta>\beta_{\rm crit}^m(\e)$.}
\subsection{The generalised eigenspace at $\lambda=0$}

In order to obtain the eigenfunction, we first differentiate \eqref{EQ:KStw2} with respect to $z$ and obtain
\begin{align*}
\begin{split}
0 &= \e u_{zzz}+cu_{zz}- w_z u^m-m w u^{m-1}u_z=\mL_{11} u_z+\mL_{12} w_z,\\
0 &= w_{zzz}+cw_{zz}-\beta(\frac{w_{z}u_z}{u})_z-\beta(\frac{wu_{zz}}{u})_z+(\beta\frac{wu_z^2}{u^2})_z=\mL_{21} u_z+\mL_{22} w_z,
\end{split}
\end{align*}
which is equivalent to
$$0=\mathcal{L}\begin{pmatrix}
u_z(z)\\ w_z(z),
\end{pmatrix}$$

where we have omitted the $\e$ subscripts and $m$ superscripts. Hence $((u_\e^m)_z,(w_\e^m)_z)$ solves the linearised eigenvalue problem \eqref{EQ:Lop_sublinear_m} for $\lambda=0$, and $\lambda=0$ is thus an eigenvalue with associated eigenfunction $((u_\e^m)_z,(w_\e^m)_z)$. This is typical for travelling wave solutions and arises from the translation invariance of solutions in the moving frame $z$.

If we instead differentiate \eqref{EQ:KStw2} with respect to the wave speed $c$, we obtain
\begin{align*}
\begin{split}
0 &= \e u_{zzc}+u_z+cu_{zc}- w_c u^m-m w u^{m-1}u_c\\
&=\mL_{11} u_c+\mL_{12} w_c+u_z,\\
0 &= w_{zzc}+w_z+cw_{zc}-\beta\left(\frac{w_{z}u_z}{u}\right)_c-\beta\left(\frac{wu_{zz}}{u}\right)_c+\beta\left(\frac{wu_z^2}{u^2}\right)_c\\
&=\mL_{21} u_c+\mL_{22} w_c+w_z,
\end{split}
\end{align*}
which is equivalent to \begin{align*}\mathcal{L}\begin{pmatrix}
u_c(z)\\ w_c(z)
\end{pmatrix}=-\begin{pmatrix}
u_z(z)\\ w_z(z)
\end{pmatrix}.\end{align*}
Hence $((u_\e^m)_c,(w_\e^m)_c)$ is a generalised eigenfunction of \eqref{EQ:Lop_sublinear_m} for $\lambda=0$ and $\lambda=0$ has algebraic multiplicity at least two and a geometric multiplicity of at least two. It was shown in \cite{harley2015numerical} that $\lambda=0$ is a second order root of the Evans function and thus the algebraic and geometric multiplicity are each precisely two.

\begin{figure}
\centering
\begin{tikzpicture}
\node at (0,0){\includegraphics[scale=1]{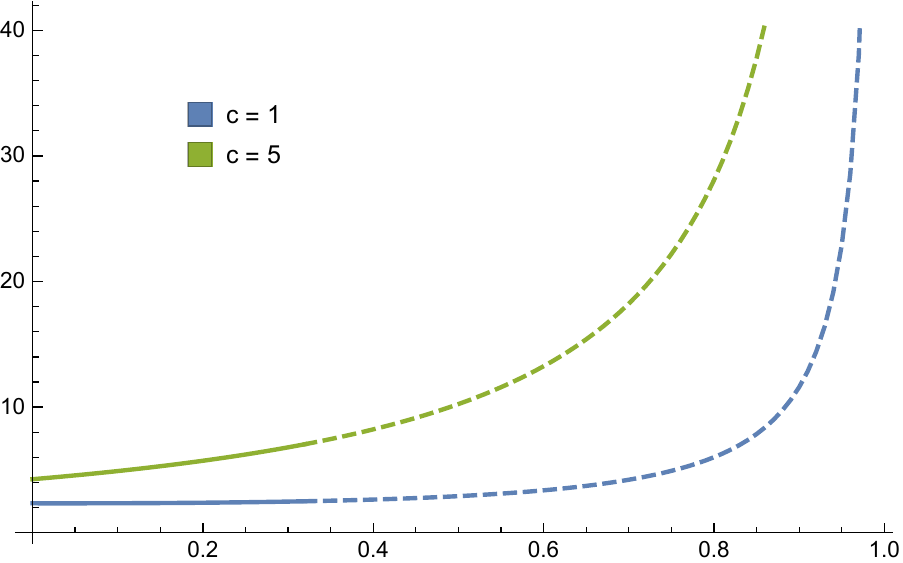}};
\node at (4.5,-3.1){$m$};
\node at (-5,0) {\rotatebox{90}{$||(u_c,w_c)||_{\mathbb{L}^2}^2$}};
\end{tikzpicture}
\caption{The $\mathbb{L}_2$ norm squared of the generalised eigenfunction for $\e=0$, \ie\ $||((u^0_m)_c,(w^0_m)_c||_{\mathbb{L}^2}^2$, is shown for $\beta=1.1$ over varying values of $m$ and two different wave speeds $c$. The solid curves represent $m$ values such that $1-m<\beta<\beta_{\rm crit}^m(0)$ and the dashed line represent values such that $\beta>\beta_{\rm crit}^m(0)$. Observe that both squared norms have an asymptote at $m=1$.}\label{FIG:norm_c}
\end{figure}

For $\e$ small we can explicitly compute the leading order (generalised) eigenfunctions from \eqref{tw_profiles_0}. In particular,
\begin{equation}\label{EQ:eigenfunctions}
\begin{split}
(u_0^m)_z(z)&=\frac{1}{c}e^{-c z} (u_0^m)^{\beta +m},\\
(w_0^m)_z(z)&=-c   e^{-cz} \left(u_0^m\right)^{\beta}+\beta e^{-cz} \left(u_0^m\right)^{\beta-1}(u_0^m)_z,\\
   (u_0^m)_c(z)&=\frac{(c z+2) u_0^m}{c \left(c^2 e^{c z}+\beta
   +m-1\right)},\\
  ( w_0^m)_c(z)&=-z   e^{-cz} \left(u_0^m\right)^{\beta}+\beta e^{-cz} \left(u_0^m\right)^{\beta-1}(u_0^m)_c.
\end{split}
\end{equation}
We have $((u_0^m)_z,(w_0^m)_z)\in \mathbb{L}^2\times \mathbb{L}^2$ and $((u_0^m)_c,(w_0^m)_c)\in\mathbb{L}^2\times \mathbb{L}^2$ with 
\begin{align}\label{EQ:uz_norm}
||((u_0^m)_z,(w_0^m)_z)||_{\mathbb{L}^2}^2=\frac{c}{4+2(\beta+m-1)}+\frac{c^5}{\beta^2-2(\beta+m-1)^2}.
\end{align}
The exact expression for $||(u_c,w_c)||_{\mathbb{L}^2}$ is not informative; instead  see Figure \ref{FIG:norm_c} for computations of $||(u_c,w_c)||_{\mathbb{L}^2}$ specific values of $\beta,m,c$. 

\subsection{Behaviour of eigenfunctions as $z\to\pm\infty$}\label{SS:boundedness}
As we have shown the existence and boundedness of the (generalised) eigenfunctions, what remains is to show that the functions are contained in the weighted spaces $\mathbb{H}^1_\nu(\mathbb{R})$ where the essential spectrum is weighted into the left half plane. We do this through examining the asymptotic behaviour of the eigenfunctions. As the exponential decay rate of the eigenfunctions as $z\to\pm\infty$ are greater than the weighted spatial eigenvalues at $z\to\pm\infty$. Hence, we conclude that $\lambda=0$ is in fact point spectrum for all $0\leq m<1$ in the space $\mathbb{H}^1_\nu(\mathbb{R})$. Since the eigenvalue problem can be expressed as an equivalent ODE \eqref{EQ:Pop_m_nonzero_eps_nonzero} in $p$, $u$ and their derivatives, it is sufficient to show $u_z$ and $u_c$ are contained in $\mathbb{H}^1_\nu(\mathbb{R})$. It was shown in \cite{wang2013mathematics} that solutions to \eqref{EQ:KStw2} can be related to solutions of a Fisher equation by first eliminating $w(z)$ using $w_\e^m(z)=e^{-cz}\left(u_\e^m(z)\right)^{\beta}$. Then, by introducing the change of variable $u_\e^m(z)=v(z)\exp\left(\frac{c}{\beta+m-1} z\right)$ the following ODE in $v(z)$ is obtained
\begin{equation}
\begin{split}\label{EQ:fisher_like}
\e v''+sv'+\eta v-v^{\beta+m}=0,\\
\end{split}
\end{equation}
where $s=c\left(1+\frac{2\e}{\beta+m-1}\right)$ and $\eta=\frac{c^2(\e+\beta+m-1)}{ (\beta+m-1)^{2}}$. We have the following result for $v(z)$ adapted from Lemma 3.1 ii a of \cite{wang2013mathematics}
\begin{lemma}(\cite{wang2013mathematics})
There exists a nonnegative travelling wave solution $v(z)$ of \eqref{EQ:fisher_like} if and only if $\beta\geq1-m$. For $\beta>1-m$ the travelling wave solution $v(z)$ is a wavefront with $v'(z)<0$ and satisfies the asymptotic conditions $$\lim_{z\to-\infty}v (z)=\eta^{\frac{1}{\beta+m-1}},\ \lim_{z\to \infty} v(z)=0.$$ The wavefront $v(z)$ has the following asymptotic behaviours: 
\begin{equation}\label{EQ:asymp_behaviour}
v(z)\sim \eta^{\frac{1}{\beta+m-1}}-C_1 e^{\kappa_1 z},\text{ as }z\to-\infty \text{ and } v(z)\sim C_2 e^{\kappa_2 z}, \text{ as }z\to\infty
\end{equation}
where
\begin{equation*}
\kappa_1=\kappa_2 -\frac{c}{2\e}+\frac{c}{2\e}\sqrt{1+\frac{4 \e (\beta +m) (\e+\beta +m-1)}{(\beta +m-1)^2}}\text{ and } 
\kappa_2 =-\frac{c}{\beta+m-1}.
\end{equation*}
\end{lemma}
From \eqref{EQ:characteristic_eqn} the unweighted spatial eigenvalues $\mu^-$ for $\lambda=0$ are given by
\begin{align*}
\mu^-_1&=\frac{c}{\beta +m-1},\ \mu^-_2=\mu^-_1m,\\
 \mu^-_{3,4}&=-\frac{c}{2
   \e }\mp\frac{c}{2\e }  \sqrt{1+\frac{4 \e  (\beta +m) (\beta +m+\e -1)}{(\beta +m-1)^2}}.
\end{align*}
Using \eqref{EQ:asymp_behaviour} and $u_\e^m(z)=v(z)e^{\frac{c}{\beta+m-1} z}$ we have the following asymptotic behaviour as $z\to-\infty$

\begin{align*}
(u_\e^m)_z(z)&\sim \frac{c\eta^{\frac{1}{\beta+m-1}}}{\beta+m-1}e^{\frac{c}{\beta+m-1}z}-C_1\left(\kappa_1+\frac{c}{\beta+m-1}\right)e^{\left(\kappa_1+\frac{c}{\beta+m-1}\right)z}\\
&\sim \frac{c\eta^{\frac{1}{\beta+m-1}}}{\beta+m-1}e^{\mu_1^-z}-C_1\left(\kappa_1+\frac{c}{\beta+m-1}\right)e^{\mu_4^-z},\\
(u_\e^m)_c(z)&\sim\left( \frac{\eta ^{\frac{1}{\beta +m-1}}+\eta _c
   \eta ^{\frac{2-\beta-m}{\beta +m-1}}}{\beta+m-1} \right)e^{\frac{c}{\beta+m-1}z}-C_1 \left(\frac{1}{\beta +m-1}+(\kappa_1)_c\right)z e^{\left(\kappa_1+\frac{c}{\beta+m-1}\right) z}\\
&\sim \left( \frac{\eta ^{\frac{1}{\beta +m-1}}+\eta _c
   \eta ^{\frac{2-\beta-m}{\beta +m-1}}}{\beta+m-1} \right)e^{\mu_1^-z} -C_1 \left(\frac{1}{\beta +m-1}+(\kappa_1)_c\right)z e^{\mu_4^- z}.
\end{align*}
The asymptotic decay rates of $u_z,$ $u_c$ as $z\to-\infty$ are precisely the two largest unstable unweighted spatial eigenvalues. In the weighted space $\mathbb{H}^1_\nu(\mathbb{R})$ the spatial eigenvalues are $\mu_-^i+\nu_-$ for $i=1,2,3,4$. It was shown in \cite{davis2017absolute} that the range of weights for $z\to-\infty$ are negative. Thus the eigenfunctions decay faster than the weighted spatial eigenvalues as $z\to-\infty$. 

From \eqref{EQ:characteristic_eqn} we have the unweighted spatial eigenvalues as $z\to\infty$ for $\lambda=0$,
\begin{align*}
&\mu_+^1=-c, && \mu_+^2=-\frac{c}{\e}, &&\mu_+^{3,4}=0.
\end{align*}
From \eqref{EQ:asymp_behaviour} we have $u_\e^m(z)\sim C_2 \text{, as }z\to\infty$ and so we can not compare the exact asymptotic exponential decay rates of the derivatives of $u_\e^m(z)$ to the spatial eigenvalues. It was shown in \cite{davis2017absolute} that positive weights $\nu_+\in(0,c)$ can be used to weight the essential spectrum into the open left half plane. Thus, as $(u_\e^m)_z$ and $(u_\e^m)_c$ are solutions to the eigenvalue problem both (generalised) eigenfunctions will decay to zero exponentially in the stable subspace spanned by the eigenvectors associated with $\mu_+^{1,2}$. Thus as $\nu_+>0$ we can conclude that there is some weighted space such that the (generalised) eigenfunctions decay faster than the weighted spatial eigenvalues as $z\to\infty$.

Hence, as the (generalised) eigenfunctions decay faster than the weighted spatial eigenvalues as $z\to\pm \infty$  we can conclude that $(u_c,w_c),(u_z,w_z)\in\mathbb{H}^1_\nu(\mathbb{R})$ for the range of weights that shift the essential spectrum into the open left half plane. Thus, the eigenvalue $\lambda=0$ is isolated and in the point spectrum in these weighted spaces.

The inclusion of a small diffusion parameter amounts to a perturbation of the operator and there are only a few possible ways new point eigenvalues can appear. These point eigenvalues emerge as perturbations of the eigenvalues in the $\e=0$ case or as new eigenvalues emerging from, and to leading order given by, the branch points of the absolute spectrum. 
From \cite{kapitula2013spectral} it follows that eigenvalues in the point spectrum are, to leading order, given by those in the $\e=0$ case. It has been shown via a numerical Evans function computation that there are no point eigenvalues in the open right half plane (excluding $\lambda=0$) when $\e=0$ up to $|\lambda|\sim \mathcal{O}(10^9)$ \cite{harley2015numerical}. As a result there can be no point spectra in this region for $0<\e\ll1$. The ODE \eqref{EQ:Pop_m_nonzero_eps_nonzero} varies smoothly in $\lambda$ as $m\to0$ and so the results of \cite{harley2015numerical} hold for $0<m<1$.
Thus, the only eigenvalues that can potentially destabilise the $\e=0$ spectrum are those that emerge from the branch points of the absolute spectrum. As the operator $\mL-\lambda I$ \eqref{EQ:Lop_sublinear_m} varies smoothly in $\lambda$ near $\e=0$ any eigenvalues that emerge from the branch points will be of the form
$$\lambda=\lambda_{\rm br}+C \e^2+\mathcal{O}(\e^3),$$
for some $C\in \mathbb{C}$ \cite{kapitula2013spectral}. Therefore, choosing $\beta<\beta_{\rm crit}^m(\e)$ such that $|\Re(\lambda_{\rm br})|$ is not $\mathcal{O}(\e^2)$ will prevent any emerging point spectrum from destabilising the spectrum.

\subsection{The limit $m\to1$}\label{SS:linear_m}

In the case of linear consumption, the travelling wave solutions are a pair of travelling wavefronts given, to leading order, by \eqref{tw_profiles_0} with $m=1$. These wavefronts satisfy \eqref{EQ:KStw2} and asymptote to 
\begin{equation}\label{EQ:limits2}
\lim_{z\to-\infty}(u^1_\e(z),w^1_\e(z))=\left(0,\frac{c^2}{\beta}+\e\frac{c^2}{\beta^2}\right) ,\quad\lim_{z\to\infty}(u^1_\e(z),w^1_\e(z))=(u_r,0),\end{equation} with $u_r$ scaled to one without loss of generality. Furthermore, the essential spectrum of $\mathcal{L}_\e^1$ \eqref{EQ:Lop_sublinear_m} includes the origin for all parameter values and all possible weights \eqref{EQ:two_sided_weight}, see Theorem 2.6 of \cite{davis2017absolute}. Therefore, this case is markedly different from $0\leq m<1$ and must be treated separately.

Similar to the previous analysis in \S\ref{SS:point_spectrum}, we compute the functions $((u^1_\e)_z,(w^1_\e)_z)$ and $((u^1_\e)_c,(w^1_\e)_c)$. These functions are given, to leading order, by \eqref{EQ:eigenfunctions} with $m=1$. While the function $((u^1_\e)_z,(w^1_\e)_z)$ persists as a solution to the eigenvalue problem \eqref{EQ:Lop_sublinear_m} with a finite norm $||((u^1_\e)_z,(w^1_\e)_z)||_{\mathbb{L}^2}$ given, to leading order, by \eqref{EQ:uz_norm}, the leading order function $((u^1_0)_c,(w^1_0)_c)$ is unbounded and hence $((u^1_\e)_c,(w^1_\e)_c)$ is not a solution to the generalised eigenvalue problem, see Figure \ref{FIG:norm_c}. The eigenvalue $\lambda=0$ is order one, in the sense that one of the eigenfunctions associated with $\lambda=0$ persists when $m\to1$. As there is no exponentially weighted function space such that $\lambda=0$ is isolated, it is not considered point spectrum.

 The intuitive reason for the reduction of order as $m\to1$ is that there is no longer a family of solutions in $c$, since, when $m=1$, the end state of $w_\e^1$ as $z\to-\infty$ depends on $c$, see \eqref{EQ:limits2}. Thus for fixed end states a travelling wave solution exists with a unique wave speed $c$, whereas for $0\leq m<1$ travelling wave solutions exist for any wave speed $c$. Alternatively, the reduction of order can be seen by examining the deformation of the absolute spectrum as $m\to1$. The absolute spectrum does not contain $\lambda=0$ for $0\leq m<1$ but, as $m\to 1$ the branches of absolute spectrum approach $\lambda=0$. Thus the order is reduced as $m\to1$ due to an eigenvalue disappearing into the absolute spectrum. 
\section{Discussion and Conclusion}
Barring the existence of extremely large values of $|\lambda|$ in the right half plane, we have now confirmed that the travelling wave solutions to \eqref{EQ:KSmain} are spectrally stable in an appropriately weighted function space for $1-m<\beta<\beta_{\rm crit}^m(\e)$ for $0\leq m<1$ and $0\leq \e\ll 1$. In particular, the point eigenvalue $\lambda=0$, proven to be of order two when $\e=m=0$ in \cite{harley2015numerical}, does not perturb in the $\e\neq 0$ case. This is because the existence of a continuous family of solutions in $c$ and the existence of a family of solutions due to translation invariance is preserved when $\e\neq0$.

For sectorial semilinear, parabolic operators with a spectral gap spectral stability implies the nonlinear (orbital) stability of the travelling wave solutions \cite{henry1981geometric}. Though we have concluded spectral stability for solutions to \eqref{EQ:KSmain} for $1-m<\beta<\beta_{\rm crit}^m(\e)$ we do not conclude nonlinear stability as the operator $\mL^m_\e$ is only quasilinear. In \cite{meyries2014quasi} it is shown that one can still conclude nonlinear stability for a large class of quasilinear parabolic operators under certain conditions. The Keller-Segel model studied in this manuscript does not fulfil these conditions, though potentially the analysis of \cite{meyries2014quasi} could be extended to encompass this model.

In the case that $m=1$ one can use a Hopf-Cole transformation in conjunction with energy estimates in order to prove the nonlinear (orbital) stability of travelling waves for $0\leq\e\ll1$ \cite{meyries2011local,wang2010chemotaxis}. These energy estimates have the potential to provide a bound on large $\lambda$. However, these energy estimates are notoriously difficult for specific linearised operators and the computation is further complicated by both the non-self-adjointness of the operator $\mL_\e^m$ given in \eqref{EQ:Lop_sublinear_m} and the fact that the Hopf-Cole transformation is not applicable when $0\leq m<1$.

As the absolute spectrum crosses into the right half plane away from the real line as $\beta\to\beta_{\rm crit}$ one expects oscillations near this critical parameter. As there are no instabilities arising from the point spectrum (barring extremely large values of $|\lambda|$) it is of interest to examine the dynamical implications and nature of this bifurcation as it is atypical for the absolute spectrum to cross away from the real line. Future work will examine this bifurcation, both analytically and numerically, to determine the impact on the wave speed, wavefront selection and whether oscillatory behaviour is observed. 

\section*{Acknowledgments}
PD and
PvH acknowledge support under the Australian Research Council's Discovery Early Career Researcher Award funding scheme DE140100741. PD and RM would like to thank Assoc. Prof. Zhi-An Wang for several informative discussions, as well as Assoc. Prof. Peter Kim and the University of Sydney for hosting these discussions.
\bibliographystyle{plain}
\bibliography{Spectral_Stability_Keller_Segel_Preprint}
\end{document}